\documentclass[12pt]{amsart}
\usepackage{amsfonts}
\usepackage{amsfonts,latexsym,rawfonts,amsmath,amssymb,amsthm}
\usepackage[plainpages=false]{hyperref}

\usepackage{graphicx}

\RequirePackage{color}

\numberwithin{equation}{section}

\newcommand{\beq}{\begin{equation}}
\newcommand{\eeq}{\end{equation}}
\newcommand{\beqs}{\begin{eqnarray*}}
\newcommand{\eeqs}{\end{eqnarray*}}
\newcommand{\beqn}{\begin{eqnarray}}
\newcommand{\eeqn}{\end{eqnarray}}
\newcommand{\beqa}{\begin{array}}
\newcommand{\eeqa}{\end{array}}

\def\p{\partial }

\def\Om{\Omega}
\def\pom{\p  \Omega}

\newtheorem{Proposition}{Proposition}[section]
\newtheorem{Theorem}[Proposition]{Theorem}
\newtheorem{Lemma}[Proposition]{Lemma}

\newtheorem{Corollary}[Proposition]{Corollary}

\title  { The anisotropic convexity of domains and the boundary estimate for two Monge-Amp\`ere   Equations   }

\begin{document}

\address{Ruosi Chen: Department of Mathematical sciences, Tsinghua University, Beijing 100084, China.}

\address{Huaiyu Jian: Department of Mathematical sciences, Tsinghua University, Beijing 100084, China.}

\email{
crs22@mails.tsinghua.edu.cn\; \;\; hjian@tsinghua.edu.cn   }

\thanks{This work was supported by NSFC 12141103  }


\bibliographystyle{plain}

\maketitle

\baselineskip=15.8pt
\parskip=3pt

\centerline {\bf  Ruosi Chen\ \ \ \ Huaiyu Jian }

\centerline {   Department of Mathematical sciences, Tsinghua University, Beijing 100084, China}

\vskip20pt

\noindent {\bf Abstract}:
We study the exact effect of the anisotropic convexity of domains on  the boundary estimate for two Monge-Amp\`ere   Equations: one is singular which
is from the proper affine hyperspheres with constant mean curvature; the other is  degenerate which is from the Monge-Amp\`ere  eigenvalue problem. As a result,
 we obtain the sharp boundary boundary estimates and the optimal global H\"older regularity for the two equations.


 \vskip20pt


\vskip20pt

\noindent {\bf  Running head:}  anisotropic convexity of domain and the boundary estimate for  Monge-Amp\`ere   Equation

\vskip20pt

\baselineskip=15.8pt
\parskip=3pt

\newpage

\maketitle

\baselineskip=15.8pt
\parskip=3.0pt



\section {Introduction}
In recent papers \cite{[JLL],[JW1]}, the second author with his collaborators introduced the
concept of the anisotropic convexity so as to describe the local convexity of domains exactly. Using the exact convexity they
obtained the sharp regularity of a kind of singular-degenerate Monge-Amp\`ere   Equations   $\det D^2 u= F(x, u, \nabla u)$
with $F$ like $A d_{x}^{\beta}|u|^{-\alpha}(1+|\nabla u|^{2})^{\frac{\gamma}{2}}$ for some positive constants $A, \beta$, $\alpha$ and $\gamma$. Here and below, we
always use the  $d_{x}=dist(x, \pom)$ to denote   the distance from the point $ x$ to the boundary of   domain $\Omega$.  Although this kind of equations
includes a few geometric models such as the affine hyperbolic sphere, the prescribed Gauss curvature problem, the $L_p$-Minkowski problem and the logarithmic Minkowski problem, it does exclude two important
equations: one is reduced by finding the proper affine hypersphere with constant mean curvature which is asymptotic to a cone in $R^{n+k+1}$ \cite{[CH]},
\begin{equation}
  \left\{ \
  \begin{aligned}
  \det D^2 u & = |u|^{-n-k-2}(x\cdot Du - u)^{-k}   \quad &in \quad &\Omega, \\
          u & = 0  \quad &on \quad &\partial\Omega;
\end{aligned}
\right.
\label{eq:monge1}
\end{equation}
the other is  the following degenerate   Monge-Amp\`ere   problem \cite{[Le1],[Le], [Ts]},
  \begin{equation}
   \left\{
    \begin{aligned}
    \det D^2 u & = |u|^q \quad &in \quad &\Omega, \\
            u & = 0  \quad &on \quad &\partial\Omega.
  \end{aligned}
  \right.
  \label{eq:monge2-q}
\end{equation}
When $q=n$, it is the well known Monge-Amp\`ere  eigenvalue problem \cite{[HH], [Lion], [LS],[Ts]}.
 Here $\Omega$ is a bounded convex domain in $\mathbb{R}^n (n\geq 2)$, which is assumed to contain the origin for problem \eqref{eq:monge1} such that
 $x\cdot Du - u>0$ in $\Omega$,   $k$  and $q$ are given positive constants.  Chen and Huang   \cite{[CH]} proved the existence and uniqueness of convex solution $u\in C^{\infty}(\Omega) \cap C^{\frac{1}{n+k+1}}(\overline{\Omega})$ to \eqref{eq:monge1}, while Le \cite{[Le]} improved the H\"older exponent to $\frac{k+2}{2(n+k+1)}$.

 Le \cite{[Le1]} proved  the existence and uniqueness of convex solution $u\in C^{\infty}(\Omega) \cap C(\overline{\Omega})$ to \eqref{eq:monge2-q}. Furthermore, in \cite{[Le]} he  obtained   the
 global  optimal    regularity
 for general  bounded convex domain:
 the solution to  \eqref{eq:monge2-q}
 $u\in   C^{\alpha}(\overline{\Omega})$ for each $\alpha\in (0, 1)$ when $q\geq n-2$ ( which should be optimal),  and   $u\in   C^{\beta}(\overline{\Omega})$ for each $\beta\in (0, \frac{2}{n-q})$ when $0<q< n-2$.

  Also, we should recall Caffarelli's result in the case $q=0$.  In \cite{[C]} he
  proved that if $\det D^2u \leq 1 $ in $\Om$ and $u=0$ on $\partial \Om$, then
$|u(x)|\leq C(n, \alpha, diam (\Om ))d_x^\alpha$ for each $ \alpha \in (0, 1)$ when $n=2$ and for $\alpha=\frac{2}{n}$ when $n\geq 3$. This result will be extended in our Theorem 1.4 (1), and it motivates
the statement of our Theorems 1.1 and 1.3 on sub-solutions and super-solutions.

However, all the above boundary estimates for the solutions are the same whether the  boundary point is a flat point or very convex point like a cone vertex.

This note is concerned with exact boundary estimate for  problems  \eqref{eq:monge1} and  for \eqref{eq:monge2-q} in the case $0\leq q< n$. We will find the exact effect of the    boundary point convexity
on the boundary estimate    for  the convex viscosity or Aleksandrov solutions.  In particular, we will obtain the global sharp  regularity  for  problems  \eqref{eq:monge1} and  \eqref{eq:monge2-q}. For this purpose, we need the following
definitions introduced in \cite{[JLL],[JW1]}.

{\bf Definition 1.1}
Suppose $\Omega$ is a bounded convex domain in $\mathbb{R}^n, x_0 \in \partial \Omega$, and $a_i \geq 1$ for $i=1,2, \cdots, n-1$. We say that {\sl $\Omega$ satisfies $\left(a_1, \cdots, a_{n-1}\right)$-exterior (-interior, resp.) convex condition at $x_0$}  if there exist positive constants $\eta_1, \cdots, \eta_{n-1}$ ( $\eta_1, \cdots, \eta_{n-1}, h,$ resp.), after suitable translation and rotation transforms, such that
$$ x_0=\boldsymbol{0} \text { and } \Omega \subseteq\left\{\left.x \in \mathbb{R}^n\left|x_n>\eta_1\right| x_1\right|^{a_1}+\cdots+\eta_{n-1}\left|x_{n-1}\right|^{a_{n-1}}\right\} .
$$
$$(
x_0= \boldsymbol{0} \text { and }\left\{x \in \mathbb{R}^n|\, \eta_1 |x_1|^{a_1}+ \cdots+\eta_{n-1}|x_{n-1}|^{a_{n-1}}<x_n<h\right\} \subseteq \Omega \subseteq \mathbb{R}_{+}^n , \; \text{resp.} )$$
If $\Omega$ satisfies $\left(a_1, \cdots, a_{n-1}\right)$-exterior (-interior, resp.) condition at every $x_0 \in \partial \Omega$ with the same positive constants $\eta_1, \cdots, \eta_{n-1}$
($\eta_1, \cdots, \eta_{n-1}, h$, resp.), we say that {\sl $\Omega$ is $\left(a_1, \cdots, a_{n-1}\right)$-exterior (-interior, resp.) convex domain.}

The convexity of $\Omega$ guarantees that the numbers $a_1, \cdots, a_{n-1}$ should be no less than 1. The less is $a_i$, the more convex is $\Omega$ along the $x_i$-direction. When $a_i$ tends to $+\infty$, $\Omega$
becomes  flat near $x_0$ along the $x_i$-direction, since $\eta_i$ may be small enough such that  $|x|<\frac{1}{2}(\eta_i)^{-\frac{1}{a_i}} $ for all $x\in \Om$, and
 the term $ \eta_i |x_i|^{a_i}$ tends to $0$. It is easy to verify that $(2, \cdots, 2)$-exterior  (-interior, resp.) condition is equivalent to the usual exterior  (-interior, resp.) sphere condition. See \cite{[JL]} for the details.

The main results of  this paper are the following four theorems.  We refer to Section 2 for the  convex viscosity (super-, sub-) solution  of problem (1.1) and Aleksandrov (super-, sub-)solution of  problem (1.2).

\begin{Theorem}
Suppose that $\Omega$ is a bounded convex domain in $\mathbb{R}^n$, $\boldsymbol{0}\in \Om$, $x_0 \in \partial\Omega$, and $\Omega$ satisfies $\left(a_1, \cdots, a_{n-1}\right)$-exterior convex condition at $x_0 $ with $a_1,\cdots,a_{n-1}\geq 2 $. If $u\in C(\overline{\Omega})$ is a   convex viscosity super-solution to problem \eqref{eq:monge1}, then there exists a positive constant $C$, depending only on $ n,k,\operatorname{diam}\Omega$ and $a_1, \cdots, a_{n-1}, \eta_1, \cdots, \eta_{n-1}$ in the exterior convex condition at $x_0$, such that
\begin{equation}
    |u(x)| \leq C d_x^{\,\theta} \text{ for all } x \in \Omega \text{ such that } d_x=\left|x-x_0\right|,
    \label{neq:theorem1}
\end{equation}
where
\begin{equation}
    \theta=\frac{\bar a+2+k}{2n+2k+2}, \ \  \bar a := \sum_{i=1}^{n-1} \frac{2}{a_i} .
    \label{eq:theta-def1}
\end{equation}
\label{theorem1}
\end{Theorem}

The following global H\"{o}lder regularity follows directly from Theorem \ref{theorem1} and Lemma 2.3 in \cite{[JL]}.

 \begin{Corollary} \label {1.2}
     Suppose that $\boldsymbol{0}\in \Om$ and  $\Omega$ is $(a_1, \cdots, a_{n-1})$-exterior convex domain in $\mathbb{R}^n$ with $a_1,\cdots,a_{n-1}\geq 2 $. If $u\in  C(\overline{\Omega})$ is a  convex viscosity super-solution to problem \eqref{eq:monge1}, then $u \in C^\theta(\bar{\Omega})$ and there exists a positive constant $C=C( n,k,\Omega)$ such that $|u|_{C^\theta(\bar{\Omega})} \leq C$, where $\theta$ is given by \eqref{eq:theta-def1}.
\end{Corollary}

\begin{Theorem}
    Suppose that $\Omega$ is a bounded convex domain in $\mathbb{R}^n,$ $\boldsymbol{0}\in \Om$, $ x_0 \in \partial\Omega$, $\Omega$ satisfies $\left(a_1, \cdots, a_{n-1}\right)$-interior convex condition at $x_0$,
  and  $\theta$  is given by \eqref{eq:theta-def1}. If $u\in  C(\overline{\Omega})$
    is convex viscosity  sub-solution to problem \eqref{eq:monge1}, then
    \begin{equation}
        |u(x)| \geq C \, d_x ^{\,\theta} \text{ for all } x \in \Omega \text{ such that } d_x=\left|x-x_0\right|< \frac{3}{8}h
        \label{neq:theorem2}
    \end{equation}
    for some positive constant $C$, depending only on $n,k,\operatorname{diam}\Omega$ and $a_1, \cdots, a_{n-1}, \eta_1, \cdots, \eta_{n-1},h$ in the interior convex condition at $x_0$.
    \label{theorem2}
\end{Theorem}
As a direct consequence of Theorem \ref{theorem2}, we have
\begin{equation} 
    |u(x)| \geq C \, d_x ^{\,\theta} \text{ for all } x \in \Omega \text{ such that } d_x< \frac{3}{8} h
\end{equation}
if  $\boldsymbol{0}\in \Om$ and  $\Omega$ is $\left(a_1, \cdots, a_{n-1}\right)$-interior convex domain.

\begin{Theorem} 
    \label{theorem3} Suppose that $\Omega$ is a bounded convex domain in $\mathbb{R}^n$,  $q\in [0, n)$, $a_i\geq 1$
    for $i=1, \cdots, n-1$,   and denote
        \begin{equation}
        \label{def-bar-a-alpha}
        \alpha = \frac{\bar{a}+2}{n-q}, \; \text{where}\;
        \bar{a} = \sum_{i=1}^{n-1} \frac{2}{a_i}\;\; \text{ as in } (1.4).
    \end{equation}
\begin{itemize}
\item[\bf (1)]
       If  $u\in  C(\overline{\Omega})$
    is  Aleksandrov  super-solution to problem \eqref{eq:monge2-q}, $ x_0 \in \partial\Omega$, and $\Omega$ satisfies $(a_1, \cdots, a_{n-1})$-exterior  condition at $x_0$, then
     \begin{equation} 
    |u(x)| \leq C d_x^{\,\lambda} \text{ for all } x \in \Omega \text{ such that } d_x=\left|x-x_0\right|,
   \end{equation} where $\lambda=\frac{\bar a+2}{n}$ for the case  $n-\bar a-2>0$,   $\lambda $ can be taken any number in $(0, 1)$ for the case $ n-\bar a -2\leq 0$,
     and $C$ is a positive constant, depending only on $n,q,\operatorname{diam}\Omega,  \lambda$ and $a_1, \cdots, a_{n-1}, \eta_1, \cdots, \eta_{n-1}$ in the exterior convex condition at $x_0$.
 \item [\bf (2)]
         If  $u\in  C(\overline{\Omega})$
    is  Aleksandrov  solution to problem \eqref{eq:monge2-q},  and $\Omega$ is $(a_1, \cdots, a_{n-1})$-exterior  convex domain, then
     \begin{equation} 
    |u(x)| \leq C d_x^{\,\lambda},  \ \ \;  \forall x \in \Omega,
    \label{neq:theorem3}
   \end{equation} where $\lambda$ can be taken any number in $(0, \alpha)$  for the case $q\in (0, n-\bar a -2)$,  $\lambda $ can be taken any number in $(0, 1)$ for the case $q\geq  \max\{n-\bar a -2, \;  0\}$,
    and $C$ is a positive constant, depending only on $n,q,\operatorname{diam}\Omega, \lambda$ and $a_1, \cdots, a_{n-1}, \eta_1, \cdots, \eta_{n-1}$ in the exterior convex condition.
    \end{itemize}
        \end{Theorem}

\begin{Theorem} 

         Suppose that  $\Omega$ is $(a_1, \cdots, a_{n-1})$-interior  convex domain,    $a_i, \bar a $ and $\alpha $ are the same as in Theorem 1.4,
         and $u\in  C(\overline{\Omega})$
    is  Aleksandrov  solution to problem \eqref{eq:monge2-q}. If
         $q\in (0, n-\bar a-2) $,  then
     \begin{equation} 
    |u(x)| \geq C d_x^{\,\lambda},  \ \ \; \forall x \in \Omega \text{ such that } d_x<\frac{3}{8}h,
    \label{neq:theorem4}
   \end{equation} where $\lambda$ can be taken any number in $(\alpha, 1)$
    and $C$ is a positive constant, depending on $n,q,\operatorname{diam}\Omega, \lambda$ and $a_1, \cdots, a_{n-1}, \eta_1, \cdots, \eta_{n-1}, h$ in the interior convex condition.
  \label{theorem4}
    \end{Theorem}

    The rest of this paper is organized as follows. In Section 2, we  recall the concept of convex viscosity solution and Aleksandrov solution to   Monge-Amp\`ere   equations, and define
 convex viscosity super-solution and sub-solution  of problem (1.1) and Aleksandrov super-solution and sub-solution of  problem (1.2).
In Sections 3, 4, 5, and 6,  we  prove Theorems 1.1, 1.3, 1.4 and 1.5 respectively.    As we see, the argument to construct smooth  super-solutions and sub-solutions based on the  anisotropic convexity of domains
 will be  delicate and technical.

\section {  viscosity solutions and Aleksandrov solutions}

  We state the definition for convex viscosity solutions to problem (1.1) and Aleksandrov solutions to problem (1.2),   which
should be well known to specialists. One can see \cite{[Cr]} and the references therein  on viscous solution theory for general elliptic equations, and \cite{[F], [G]}
on Aleksandrov solutions.

{\bf Definition 2.1}\; Let $\Omega \subset \mathbb{R}^n$ be a convex domain, $\boldsymbol{0}\in \Om$,        $u\in C(\overline{\Omega})$ be convex in  $\Omega$,
and $u<0$ in $\Om$.
{\sl We say that   $u$ is a  convex  viscosity sub-solution (super-solution) of the equation
\begin{equation}\label{2.1}\det D^2u=|u|^{-n-k-2}(x\cdot Du - u)^{-k}\;\;   in \;\; \Omega
 \end{equation} if whenever convex function $\phi \in C^2(\Omega) $
and $x_0\in \Omega$ are such that $(u-\phi)(x)\leq (\geq ) (u-\phi)(x_0)$ for all $x$ in a neighborhood of $x_0$,
 then we must have
$$\det D^2\phi (x_0)\geq (\leq ) |u (x_0)|^{-n-k-2}(x_0\cdot D\phi(x_0) - u(x_0))^{-k}.$$ }   {\sl If $u$ is both convex viscosity sub-solution and viscosity super-solution, then $u$ is called a viscosity solution.}

{\bf Definition 2.2}\; {\sl A convex viscosity sub-solution (super-solution) $u$ to (2.1) is called   convex viscosity  sub-solution (super-solution)  to problem (1.1) if $u=0$ on $\partial \Om$.}

Denote $$
	\partial u(x_0)=\left\{p \in \mathbb{R}^{n} |\; u(x) \geq u(x_0)+p \cdot(x-x_0) \text{ for all } x \in \Omega\right\}.
	$$   Obviously, $u(x_0)+p \cdot(x-x_0) $ is a supporting hyperplane function of $u$ at $x_0$ if $p \in \partial u(x_0)$.    Let $\partial u(E):=\bigcup_{x \in E} \partial u(x)$. It is measurable for any open and closed subset $E \subset \Omega$, and so for any Borel set $E \subset  \Omega$.

{\bf Definition 2.3}\;    Let $u\in C(\overline{\Omega})$ be convex in  $\Omega$. {\sl The  Monge-Amp\`ere measure $\mathbb{M} u$ is defined by $\mathbb{M} u(E)=|\partial u(E)|$ for each Borel set $E \subset \Omega$, $u$ is called  Aleksandrov sub-solution (super-solution, solution) to
   \begin{equation}\label{2.2}\det D^2 u = |u|^q \ \ \text{  in}\ \  \Om
     \end{equation} if $\mathbb{M} u\geq (\leq , \; =) |u(x)|^qdx$  in $\Om$.}

     {\bf Definition 2.4}\; {\sl An Aleksandrov sub-solution (super-solution, solution) $u$ to (2.2) is called  Aleksandrov sub-solution (super-solution, solution)  to problem (1.2) if $u=0$ on $\partial \Om$.}


It is obvious that a convex $C^2$-solution    must be   convex viscosity solution and Aleksandrov solution, and $\mathbb{M} u=\det D^2u$.   Caffarelli   indicated that the Aleksandrov solution   to the   equation  $\det D^2u=\eta(x)$ is equivalent to the viscosity solution if $\eta\in C(\Omega)$. See \cite {[G]} for the proof.
	
	We will need the comparison principle for Aleksandrov  solutions. See,  for example, \cite{[F], [G]}  for the details.
	
	\begin{Lemma} 
		Suppose that $\Om $ is a bounded convex domain, $u,v \in C(\overline{\Om})$, $u \geq v$   on  $\partial \Om$    and
		$$
	   \mathbb{M} u(E) \leq  \mathbb{M}v(E) < \infty \text{ for all Borel set }  E \subset \Om.
		$$
		Then $
		u(x) \geq v(x) $    in  $\Om$.
		
	\end{Lemma}

\section{The proof of Theorem \ref{theorem1} \label{pfthm1}}

Throughout this section, we suppose the assumption of Theorem 1.1.  We start proving this theorem.

 For  $ x_0  \in \partial\Omega $, denote $ x_0 = (y_1,\cdots,y_n)$.
    Notice that \eqref{eq:monge1} is invariant under rotation centered at the origin.  By the assumption we have
    \begin{equation*}
        \Omega \subseteq \{x\in\mathbb{R}^n \mid x_n-y_n \geq \eta_1 |x_1-y_1|^{a_1}+\cdots+\eta_{n-1} |x_{n-1}-y_{n-1}|^{a_{n-1}} \}.
    \end{equation*}
    Since $\boldsymbol{0} \in \mathbb{R}^n$ is an inner point of $\Omega$, we have
    \begin{equation}
        \label{yn-d0}
        -y_n \geq d_0 := \operatorname{dist}(0,\partial\Omega).
    \end{equation}
    Let $\widetilde{x}=x-x_0$, $v(\widetilde{x})=u(x)$ and $\widetilde{\Omega} = \Omega-x_0.$ Then  $v(x)$ is convex super-solution to the problem
    \begin{equation}
        \left\{ \
          \begin{aligned}
          \det D^2 v & = |v|^{-n-k-2}((x+x_0)\cdot Dv - v)^{-k}   \quad &in \quad &\widetilde{\Omega}, \\
                  v & = 0  \quad &on \quad &\partial\widetilde{\Omega},
        \end{aligned}
        \right.
        \label{eq:monge-v}
    \end{equation}
    and $\widetilde{\Omega}\subseteq \{x\in\mathbb{R}^n \mid x_n \geq \eta_1 |x_1|^{a_1}+\cdots +\eta_{n-1} |x_{n-1}|^{a_{n-1}} \}.$ Then for any $z\in \widetilde{\Omega} $ such that $d_z = |z-\boldsymbol{0}|, $ we have $z \bot \{x\in \mathbb{R}^n \mid x_n = 0\}$ and thus $d_z = z_n. $  It suffices to prove
    \begin{equation}
        \label{eq:theorem1-goal}
        |v(0,\cdots,0,z_n)| \leq C z_n^{\theta}, \quad \forall z \in \widetilde{\Omega},
    \end{equation}
    which is equivalent to (1.3).

    For all $i \in\{1,2, \cdots, n-1\}$, we have
    $$
        \begin{aligned}
            \widetilde{\Omega} & \subseteq\left\{ x \in \mathbb{R}^n\,|\, x_n>\eta_1 |x_1|^{a_1}+\cdots+\eta_{n-1}\left|x_{n-1}\right|^{a_{n-1}}\right\} \\
            & \subseteq\left\{ x \in \mathbb{R}^n \, |\, x_n>\eta_i |x_i|^{a_i}\right\}=\{x \in \mathbb{R}^n \, |\, (\frac{\varepsilon}{\eta_i} )^{\frac{2}{a_i}} (\frac{x_n}{\varepsilon})^{\frac{2}{a_i}}> |x_i|^2 \},
        \end{aligned}
    $$
    where $\varepsilon \in (0, \min \limits_{1 \leq i \leq n-1} \eta_i )$ which will be specified below. Denote
    \begin{equation}
        \delta(\varepsilon)=\max _{1 \leq i \leq n-1} (\frac{\varepsilon}{\eta_i} )^{\frac{2}{a_i}},
        \label{delta-def}
    \end{equation}
    then for all $i \in\{1,2, \cdots, n-1\}$, we have
    \begin{equation}
        \widetilde{\Omega} \subseteq \{ x \in \mathbb{R}^n\,|\, \delta(\varepsilon)(\frac{x_n}{\varepsilon})^{\frac{2}{a_i}}> |x_i|^2 \} .
        \label{omega-sub}
    \end{equation}
    Again for all $i \in\{1,2, \cdots, n-1\}$, let
    \begin{equation}
        \begin{aligned}
        & W_i(x)=-\left[\left(\frac{x_n}{\varepsilon}\right)^{\frac{2}{a_i}}-x_i^2\right]^{\frac{1}{b_i}}, \quad b_i=\frac{2}{a_i\theta}, \quad \forall x \in \widetilde{\Omega}, \\
        & W(x)=\sum_{i=1}^{n-1} W_i, \quad \forall x \in \widetilde{\Omega}.
        \end{aligned}
        \label{eq:def-Wi-W}
    \end{equation}
     Then
    \begin{equation*}
        \left|W_i\right|^{b_i} =\left(\frac{x_n}{\varepsilon}\right)^{\frac{2}{a_i}}-x_i^2 \in [(1-\delta(\varepsilon)) \cdot (\frac{x_n}{\varepsilon} )^{\frac{2}{a_i}}, \, (\frac{x_n}{\varepsilon} )^{\frac{2}{a_i}}],
    \end{equation*}
    which is equivalent to
    \begin{equation}
        \left|W_i\right| \in [(1-\delta(\varepsilon))^{\frac{1}{b_i}} \cdot (\frac{x_n}{\varepsilon} )^\theta, \,(\frac{x_n}{\varepsilon} )^\theta].
        \label{Wi-range}
    \end{equation}
    It is also clear that $W\in C^2(\Om)\cap C(\overline{\Om})$ and $W\leq 0$ on $\partial\widetilde{\Omega}$.  Next, we will prove that $W$ is convex in $\Om$ and satisfies
    \begin{equation}
    \label{eq:goal}
    H[W] := \det D^2W \cdot |W|^{n+k+2}\,((x+x_0)\cdot DW-W)^k >1 \ \ in \ \ \widetilde{\Omega}.
    \end{equation}
    The  calculations from now to (3.18) is almost the same as the proof of Theorem 1.1 in \cite{[JLL]}, since the function $W$ have the same form except for the different $\theta$.
     But for convenience,  we write the details here.    First, we have
    \begin{equation}
        \begin{aligned}
        W_{x_i x_i}& =  \left(W_i\right)_{x_i x_i} = \frac{2}{b_i}\left|W_i\right|^{1-b_i}+\frac{4\left(b_i-1\right)}{b_i^2}\left|W_i\right|^{1-2 b_i} x_i^2, \quad \forall i \in\{1,2, \cdots, n-1\}, \\
        W_{x_i x_j} &= 0, \quad \forall i, j \in\{1,2, \cdots, n-1\},\quad  i \neq j \\
        W_{x_n x_n} &= \sum_{i=1}^{n-1} \left[\,\frac{2\left(a_i-2\right)}{a_i^2 b_i}\left|W_i\right|^{1-b_i}\left(\frac{x_n}{\varepsilon}\right)^{\frac{2}{a_i}-2}(\frac{1}{\varepsilon})^2  +\frac{4\left(b_i-1\right)}{a_i^2 b_i^2}\left|W_i\right|^{1-2 b_i}\left(\frac{x_n}{\varepsilon}\right)^{\frac{4}{a_i}-2}(\frac{1}{\varepsilon})^2\, \right], \\
        W_{x_i x_n} & = \left(W_i\right)_{x_i x_n}=-\frac{4\left(b_i-1\right)}{a_i b_i^2}\left|W_i\right|^{1-2 b_i}\left(\frac{x_n}{\varepsilon}\right)^{\frac{2}{a_i}-1} \cdot \frac{x_i}{\varepsilon}, \quad \forall i \in\{1,2, \cdots, n-1\} .
        \end{aligned}
        \label{dWs}
    \end{equation}
    No matter $b_i \geq 1$ or $b_i<1$, by \eqref{omega-sub},\eqref{Wi-range} and \eqref{dWs}, we always have
    \begin{equation} 
        \begin{aligned}
        W_{x_i x_i} & \geq  \frac{2}{b_i} \cdot \min \{(1-\delta(\varepsilon))^{\frac{1-b_i}{b_i}}, 1\} \cdot(\frac{x_n}{\varepsilon})^{\theta(1-b_i)} \\
        &\ \ -\frac{4|b_i-1|}{b_i^2} \cdot \max \{(1-\delta(\varepsilon))^{\frac{1-2 b_i}{b_i}}, 1\} \cdot(\frac{x_n}{\varepsilon})^{\theta(1-2 b_i)} \cdot \delta(\varepsilon) \cdot(\frac{x_n}{\varepsilon})^{\theta b_i} \\
        & =  \frac{2}{b_i}[\min \{(1-\delta(\varepsilon))^{\frac{1-b_i}{b_i}}, 1\}-\delta(\varepsilon) \frac{2|b_i-1|}{b_i} \max \{(1-\delta(\varepsilon))^{\frac{1-2 b_i}{b_i}}, 1\}] \\
        & \ \ \cdot(\frac{x_n}{\varepsilon})^{\theta(1-b_i)} \\
        &:=  c_i(\varepsilon) \cdot(\frac{x_n}{\varepsilon})^{\theta(1-b_i)} .
        \end{aligned}
        \label{Wxixi}
    \end{equation}
    Since $\lim\limits_{\varepsilon \to 0} \delta(\varepsilon) = 0 $ by \eqref{delta-def}, we have
    \begin{equation}
        \label{ci-to-0}
        \lim_{\varepsilon \to 0} c_i(\varepsilon) = \frac{2}{b_i} = a_i \theta >0.
    \end{equation}
    Denote $ \xi_i(x) = |W_i(x)| (\frac{x_n}{\varepsilon})^{-\theta} $ for all $i\in\{1,2,\cdots,n-1\}.$ It follows from \eqref{Wi-range} that
    \begin{equation*}
        \xi_i(x) \in [(1-\delta(\varepsilon))^{\frac{1}{b_i}},1], \quad \forall x \in \widetilde{\Omega}, \quad \forall i \in \{1,2,\cdots,n-1\},
    \end{equation*}
    and therefore
    \begin{equation}
        \label{xii-to-0}
        \lim_{\varepsilon \to 0} \xi_i(x) = 1
    \end{equation}
    uniformly for $x\in \widetilde{\Omega}.$ Then
    \begin{equation}
        \label{Wxnxn}
        \begin{aligned}
        W_{x_n x_n} & =\sum_{i=1}^{n-1}\,[\frac{2(a_i-2)}{a_i^2 b_i} \xi_i^{1-b_i}(\frac{x_n}{\varepsilon})^{\theta-2}(\frac{1}{\varepsilon})^2+\frac{4(b_i-1)}{a_i^2 b_i^2} \xi_i^{1-2 b_i}(\frac{x_n}{\varepsilon})^{\theta-2}(\frac{1}{\varepsilon})^2] \\
        & = (\frac{1}{\varepsilon})^2 \cdot \theta^2 \sum_{i=1}^{n-1}\,[(\frac{1}{\theta}-b_i) \xi_i^{1-b_i}+(b_i-1) \xi_i^{1-2 b_i}] \cdot(\frac{x_n}{\varepsilon})^{\theta-2} \\
        & :=(\frac{1}{\varepsilon})^2 \cdot c_n(\varepsilon) \cdot(\frac{x_n}{\varepsilon})^{\theta-2},
        \end{aligned}
    \end{equation}
    where $c_n(\varepsilon)$ satisfies
    \begin{equation}
        \label{cn-to-0}
        \lim_{\varepsilon \to 0} c_n(\varepsilon) = (n-1)(\frac{1}{\theta}-1) \theta^2 > 0.
    \end{equation}
    Similarly, we have
    \begin{equation}
        \label{Wxixn}
        \begin{aligned}
        |W_{x_i x_n}| & \leq \frac{4|b_i-1|}{a_i b_i^2} \cdot \xi_i^{1-2 b_i} \cdot(\frac{x_n}{\varepsilon})^{\theta-\frac{2}{a_i}-1} \cdot(\delta(\varepsilon))^{\frac{1}{2}} \cdot(\frac{x_n}{\varepsilon})^{\frac{1}{a_i}} \cdot \frac{1}{\varepsilon} \\
        & =(\delta(\varepsilon))^{\frac{1}{2}} \cdot \frac{1}{\varepsilon} \cdot \frac{4|b_i-1|}{a_i b_i^2} \xi_i^{1-2 b_i} \cdot(\frac{x_n}{\varepsilon})^{\theta-\frac{1}{a_i}-1} \\
        & :=(\delta(\varepsilon))^{\frac{1}{2}} \cdot \frac{1}{\varepsilon} \cdot \widetilde{c_i}(\varepsilon) \cdot(\frac{x_n}{\varepsilon})^{\theta-\frac{1}{a_i}-1},
        \end{aligned}
    \end{equation}
    and it follows from \eqref{xii-to-0} that
    \begin{equation}
        \label{ci-tilde-to-0}
        \lim_{\varepsilon \to 0} \widetilde{c_i}(\varepsilon) = \frac{4|b_i-1|}{a_i b_i^2}
    \end{equation}
    uniformly for $x \in \widetilde{\Omega}.$ From \eqref{dWs}, we can write
    $$
    D^2 W:=\left(\begin{array}{cc}
    A_{n-1} & \vec{v} \\
    \vec{v}^T & W_{x_n x_n}
    \end{array}\right),
    $$
    where $A_{n-1}=\operatorname{diag} (W_{x_1 x_1}, \cdots, W_{x_{n-1} x_{n-1}} )$ and $\vec{v}^T=(W_{x_1 x_n}, \cdots, W_{x_{n-1} x_n})$. Then
    $$
    \operatorname{det} D^2 W=\operatorname{det} A_{n-1} \cdot(W_{x_n x_n}-\vec{v}^T A_{n-1}^{-1} \vec{v}\,) .
    $$
    Using \eqref{Wxixi} and \eqref{Wxixn}, we get
    $$
    \begin{aligned}
    \vec{v}^T A_{n-1}^{-1} \vec{v} & =\sum_{i=1}^{n-1} \frac{[W_{x_i x_n}]^2}{W_{x_i x_i}} \\
    & \leq \sum_{i=1}^{n-1} \frac{[(\delta(\varepsilon))^{\frac{1}{2}} \cdot \frac{1}{\varepsilon} \cdot \widetilde{c_i}(\varepsilon) \cdot(\frac{x_n}{\varepsilon})^{\theta-\frac{1}{a_i}-1}]^2}{c_i(\varepsilon) \cdot (\frac{x_n}{\varepsilon})^{\theta(1-b_i)}} \\
    & =\delta(\varepsilon) \cdot(\frac{1}{\varepsilon})^2 \cdot \sum_{i=1}^{n-1} \frac{(\widetilde{c_i}(\varepsilon))^2}{c_i(\varepsilon)} \cdot(\frac{x_n}{\varepsilon})^{\theta-2},
    \end{aligned}
    $$
    which, combined with \eqref{Wxnxn}, implies that
    \begin{equation}
        \label{est-detD2W}
        \begin{aligned}
        \operatorname{det} D^2 W &\geq  c_1(\varepsilon) \cdot(\frac{x_n}{\varepsilon})^{\theta(1-b_1)} \cdots c_{n-1}(\varepsilon) \cdot(\frac{x_n}{\varepsilon})^{\theta(1-b_{n-1})} \\
        &\ \ \cdot(\frac{1}{\varepsilon})^2[c_n(\varepsilon) \cdot(\frac{x_n}{\varepsilon})^{\theta-2}-\delta(\varepsilon) \cdot \sum_{i=1}^{n-1} \frac{(\widetilde{c_i}(\varepsilon))^2}{c_i(\varepsilon)} \cdot(\frac{x_n}{\varepsilon})^{\theta-2}] \\
        &=  (\frac{1}{\varepsilon})^2 \cdot c_1(\varepsilon) \cdots c_{n-1}(\varepsilon) \cdot[c_n(\varepsilon)-\delta(\varepsilon) \cdot \sum_{i=1}^{n-1} \frac{(\widetilde{c_i}(\varepsilon))^2}{c_i(\varepsilon)}] \\
        & \ \  \cdot(\frac{x_n}{\varepsilon})^{\theta(n-b_1-b_2 \cdots-b_{n-1})-2} \\
        &:=  (\frac{1}{\varepsilon})^2 \cdot \tau_1(\varepsilon) \cdot(\frac{x_n}{\varepsilon})^{n \theta-(\frac{2}{a_1}+\cdots+\frac{2}{a_{n-1}})-2},
        \end{aligned}
    \end{equation}
    and it follows from \eqref{ci-to-0},\eqref{cn-to-0} and \eqref{ci-tilde-to-0} that $\tau_1(\varepsilon)$ satisfies
    $$
    \begin{aligned}
    \lim _{\varepsilon \to 0} \tau_1(\varepsilon) & =\lim _{\varepsilon \to 0} c_1(\varepsilon) \cdots c_{n-1}(\varepsilon)\left[c_n(\varepsilon)-\delta(\varepsilon) \sum_{i=1}^{n-1} \frac{(\widetilde{c_i}(\varepsilon))^2}{c_i(\varepsilon)}\right] \\
    & =\lim _{\varepsilon \to 0} c_1(\varepsilon) \cdots \lim _{\varepsilon \to 0} c_{n-1}(\varepsilon) ( \lim _{\varepsilon \to 0} c_n(\varepsilon)-0 ) \\
    & = a_1 \cdots a_{n-1} \theta^{n-1}(n-1)(\frac{1}{\theta}-1) \theta^2 \\
    & = (n-1)a_1 \cdots a_{n-1} \theta^n(1-\theta) >0 \ \ \text{ uniformly for } x \in \widetilde{\Omega}.
    \end{aligned}
    $$
    From \eqref{Wi-range}, we get $|W_i| \geq (1-\delta(\varepsilon))^{\frac{1}{b_i}} \cdot (\frac{x_n}{\varepsilon} )^\theta$, then
    $$
    |W| \geq (n-1) \min_{1\leq i \leq n-1} [(1-\delta(\varepsilon))^{\frac{1}{b_i}}] \cdot (\frac{x_n}{\varepsilon} )^\theta,
    $$
    which yields
    \begin{equation}
        \label{W^n+k+2}
        |W|^{n+k+2} \geq \left[(n-1) \min_{1\leq i \leq n-1} [(1-\delta(\varepsilon))^{\frac{1}{b_i}}]\right]^{n+k+2} \cdot (\frac{x_n}{\varepsilon} )^{\theta(n+k+2)}.
    \end{equation}

    Next we estimate $(x+x_0)\cdot DW-W .$ Using (3.6), by a direct calculation we have
    \begin{equation*}
        \begin{aligned}
          (x+x_0)&\cdot DW-W = (x_1+y_1,\cdots,x_n+y_n)\cdot DW - W \\
            =\, & \sum_{i=1}^{n-1} (x_i+y_i) W_{x_i} + (x_n+y_n) W_{x_n} - \sum_{i=1}^{n-1} W_i \\
            =\, & \sum_{i=1}^{n-1} (x_i+y_i) (W_i)_{x_i} + (x_n+y_n) \sum_{i=1}^{n-1}(W_i)_{x_n} - \sum_{i=1}^{n-1} W_i \\
            =\, & \sum_{i=1}^{n-1} (x_i+y_i)\frac{2x_i}{b_i}[(\frac{x_n}{\varepsilon})^{\frac{2}{a_i}}-x_i^2]^{\frac{1}{b_i}-1}+(x_n+y_n)\,\frac{-2x_n^{\frac{2}{a_i}-1}}{a_i b_i \varepsilon^{\frac{2}{a_i}}}\,[(\frac{x_n}{\varepsilon})^{\frac{2}{a_i}}-x_i^2]^{\frac{1}{b_i}-1} + \sum_{i=1}^{n-1} |W_i| \\
            = \, & \sum_{i=1}^{n-1}\frac{2}{b_i}\,[(\frac{x_n}{\varepsilon})^{\frac{2}{a_i}}-x_i^2]^{\frac{1}{b_i}-1}\,\left[\frac{b_i}{2}[(\frac{x_n}{\varepsilon})^{\frac{2}{a_i}}-x_i^2]-\frac{x_n^\frac{2}{a_i}}{a_i \varepsilon^{\frac{2}{a_i}}} - \frac{y_n \cdot x_n^{\frac{2}{a_i}-1}}{a_i \varepsilon^{\frac{2}{a_i}}} +x_i y_i +x_i^2 \right].
        \end{aligned}
    \end{equation*}
    Since $x=(x_1,\cdots,x_n)\in \widetilde{\Omega}$, $\boldsymbol{0}\in \partial \widetilde{\Omega}$, $\boldsymbol{0}\in \Om$ and $(y_1,\cdots,y_n) \in \partial\Omega,$ we have
     $$ |x_i| \leq \operatorname{diam}\Omega, |y_i| \leq \operatorname{diam}\Omega,\ \  i=1, \cdots, n-1.$$  Therefore for $i=1, \cdots, n-1,$
    \begin{equation*}
        \begin{aligned}
            & \frac{b_i}{2}[(\frac{x_n}{\varepsilon})^{\frac{2}{a_i}}-x_i^2]-\frac{x_n^\frac{2}{a_i}}{a_i \varepsilon^{\frac{2}{a_i}}} - \frac{y_n \cdot x_n^{\frac{2}{a_i}-1}}{a_i \varepsilon^{\frac{2}{a_i}}} +x_i y_i +x_i^2 \\
            \geq\ &  \frac{b_i}{2}[(\frac{x_n}{\varepsilon})^{\frac{2}{a_i}}-x_i^2]-\frac{1}{a_i}(\frac{x_n}{\varepsilon})^{\frac{2}{a_i}} - \frac{y_n \cdot x_n^{\frac{2}{a_i}-1}}{a_i \varepsilon^{\frac{2}{a_i}}} -(\operatorname{diam}\Omega)^2.
        \end{aligned}
        \label{eq:est-sum}
    \end{equation*}
    We would like to have
    \begin{equation}
        \label{eq:yn-diam}
        \frac{- y_n \cdot x_n^{\frac{2}{a_i}-1}}{2\,a_i\, \varepsilon^{\frac{2}{a_i}}} -(\operatorname{diam}\Omega)^2 \geq 0,
    \end{equation}
    which is equivalent to $$ \varepsilon \leq (\frac{(-y_n)\,x_n^{\frac{2}{a_i}-1}}{2 a_i (\operatorname{diam}\Omega)^2 })^{\frac{a_i}{2}}, \quad \forall i\in\{1,\cdots,n-1\} . $$
    Notice that $\frac{2}{a_i}-1 \leq 0$ because of  $a_i \geq 2$, we can take
    $$\varepsilon_1 := \min\limits_{1\leq i \leq n-1} (\frac{d_0 \, (\operatorname{diam}\Omega)^{\frac{2}{a_i}-1}}{2 a_i (\operatorname{diam}\Omega)^2 })^{\frac{a_i}{2}} = C(\Omega) >0 ,$$
    where $d_0$ is the same as in  \eqref{yn-d0}. Then for $ \varepsilon < \varepsilon_1$, inequality \eqref{eq:yn-diam} holds true. Again from \eqref{Wi-range}, we get
    \begin{equation*}
        (\frac{x_n}{\varepsilon})^{\frac{2}{a_i}} \leq \frac{1}{1-\delta} |W_i|^{b_i} = \frac{1}{1-\delta} [ (\frac{x_n}{\varepsilon})^{\frac{2}{a_i}}-x_i^2] .
    \end{equation*}
    In \eqref{delta-def}, we can choose $\varepsilon < \varepsilon_2 = c(a_1,\cdots,a_{n-1}, \eta_1,\cdots,\eta_{n-1}) $ sufficiently small such that $\delta(\varepsilon) < 1-\theta = 1- \frac{2}{a_i b_i},$ then $\frac{b_i}{2} > \frac{1}{a_i(1-\delta(\varepsilon))} .$  It follows from this and (3.19) that
    \begin{equation*}
        \begin{aligned}
            & \frac{b_i}{2}[(\frac{x_n}{\varepsilon})^{\frac{2}{a_i}}-x_i^2]-\frac{1}{a_i}(\frac{x_n}{\varepsilon})^{\frac{2}{a_i}} - \frac{y_n \cdot x_n^{\frac{2}{a_i}-1}}{a_i \varepsilon^{\frac{2}{a_i}}} -(\operatorname{diam}\Omega)^2 \\
            \geq\ & \frac{b_i}{2}[(\frac{x_n}{\varepsilon})^{\frac{2}{a_i}}-x_i^2]-\frac{1}{a_i}(\frac{x_n}{\varepsilon})^{\frac{2}{a_i}} - \frac{y_n \cdot x_n^{\frac{2}{a_i}-1}}{2\, a_i\, \varepsilon^{\frac{2}{a_i}}} \\
            \geq\ & \frac{b_i}{2}[(\frac{x_n}{\varepsilon})^{\frac{2}{a_i}}-x_i^2] - \frac{1}{a_i(1-\delta(\varepsilon))}[(\frac{x_n}{\varepsilon})^{\frac{2}{a_i}}-x_i^2] - \frac{y_n \cdot x_n^{\frac{2}{a_i}-1}}{2\, a_i\, \varepsilon^{\frac{2}{a_i}}} \\
            \geq\ & \frac{(-y_n) \cdot x_n^{\frac{2}{a_i}-1}}{2\, a_i\, \varepsilon^{\frac{2}{a_i}}}.
        \end{aligned}
    \end{equation*}
    Using \eqref{Wi-range} again and combing the above three estimates,  we have
    \begin{equation}
        \label{est-xDW-W}
        \begin{aligned}
          (x+x_0)\cdot DW-W
            \geq\ & \sum_{i=1}^{n-1}\frac{2}{b_i}\,[(\frac{x_n}{\varepsilon})^{\frac{2}{a_i}}-x_i^2]^{\frac{1}{b_i}-1} \frac{(-y_n) \cdot x_n^{\frac{2}{a_i}-1}}{2\, a_i\, \varepsilon^{\frac{2}{a_i}}} \\
            =\ & (-y_n)\cdot \sum_{i=1}^{n-1} \frac{2}{a_i b_i}\, [(\frac{x_n}{\varepsilon})^{\frac{2}{a_i}}-x_i^2]^{\frac{1}{b_i}-1} \, \frac{x_n^{\frac{2}{a_i}-1}}{2 \varepsilon^{\frac{2}{a_i}}} \\
            \geq \ & (-y_n)\, \theta \cdot \sum_{i=1}^{n-1}  [(1-\delta(\varepsilon))(\frac{x_n}{\varepsilon})^{\frac{2}{a_i}}]^{\frac{1}{b_i}-1} \, \frac{x_n^{\frac{2}{a_i}-1}}{2 \varepsilon^{\frac{2}{a_i}}} \\
            \geq \ & \frac{d_0\,\theta}{\varepsilon}\left[\sum_{i=1}^{n-1} (1-\delta(\varepsilon))^{\frac{1}{b_i}-1}\right] (\frac{x_n}{\varepsilon})^{\theta-1}.
        \end{aligned}
    \end{equation}
    Combining \eqref{est-detD2W}, \eqref{W^n+k+2} and \eqref{est-xDW-W}, we finally get
    \begin{equation*}
        \begin{aligned}
            H[W] & = \det D^2W \cdot |W|^{n+k+2}\,((x+x_0)\cdot DW-W)^k \\
            & \geq (\frac{1}{\varepsilon})^2  \tau_1(\varepsilon)\cdot \left[(n-1) \min_{1\leq i \leq n-1} [(1-\delta(\varepsilon))^{\frac{1}{b_i}}]\right]^{n+k+2} \cdot \frac{d_0^k\,\theta^k}{\varepsilon^k}\left[\sum_{i=1}^{n-1} (1-\delta(\varepsilon))^{\frac{1}{b_i}-1}\right]^k\\
            & \quad \cdot (\frac{x_n}{\varepsilon})^{n \theta-(\frac{2}{a_1}+\cdots+\frac{2}{a_{n-1}})-2} \cdot (\frac{x_n}{\varepsilon} )^{\theta(n+k+2)} \cdot (\frac{x_n}{\varepsilon})^{k\theta-k} \\
            & = (\frac{1}{\varepsilon})^2  \tau_1(\varepsilon)\cdot \left[(n-1) \min_{1\leq i \leq n-1} [(1-\delta(\varepsilon))^{\frac{1}{b_i}}]\right]^{n+k+2} \cdot \frac{d_0^k\,\theta^k}{\varepsilon^k}\left[\sum_{i=1}^{n-1} (1-\delta(\varepsilon))^{\frac{1}{b_i}-1}\right]^k \\
            & \quad \cdot (\frac{x_n}{\varepsilon})^{(2n+2k+2) \theta-(\frac{2}{a_1}+\cdots+\frac{2}{a_{n-1}})-2-k}.
        \end{aligned}
    \end{equation*}
    Putting the $\theta$ given by \eqref{eq:theta-def1} into the above estimate, we have
    \begin{equation*}
        H[W] \geq (\frac{1}{\varepsilon})^2  \tau_1(\varepsilon)\cdot \left[(n-1) \min_{1\leq i \leq n-1} [(1-\delta(\varepsilon))^{\frac{1}{b_i}}]\right]^{n+k+2} \cdot \frac{d_0^k\,\theta^k}{\varepsilon^k}\left[\sum_{i=1}^{n-1} (1-\delta(\varepsilon))^{\frac{1}{b_i}-1}\right]^k .
    \end{equation*}
    Since $\tau_1(\varepsilon)$  is of  positive lower bound uniformly for $\widetilde{ x}\in \widetilde{\Om} $ when $\varepsilon \to 0^+$, we can
    choose $\varepsilon = \varepsilon_0 = \varepsilon_0(n,k,d_0,\operatorname{diam}\Omega, a_1,\cdots,a_{n-1}, \eta_1,\cdots,\eta_{n-1}) = \varepsilon_0(n,k,\Omega) > 0$ sufficiently small,   $\varepsilon_0<\min\{\varepsilon_1,\varepsilon_2 \} < \min \limits_{1 \leq i \leq n-1} \eta_i $, such that  $H[W]> 1$. Therefore \eqref{eq:goal} has been proved.

    Now we claim that $
   0 \geq v \geq W$ on $\widetilde{\Omega}$. Otherwise,  we have some $\bar x\in  \widetilde{\Omega}$ such that  $v-W$ attains a negative minimum at $\bar x$. Since $v$ is convex super-solution
   to problem (3.2), by Definition 2.1 we have
   \begin{equation*}
        \begin{aligned}
   \det D^2W (\bar x) & \leq |v (\bar x)|^{-n-k-2}(\bar x \cdot DW(\bar x) - v(\bar x))^{-k}\\
   & < |W(\bar x)|^{-n-k-2}(\bar x \cdot DW(\bar x) - W(\bar x))^{-k}
   \end{aligned}
   \end{equation*}
    which contradicts \eqref{eq:goal}.
     As the special case of the claim, $\forall z = (0,\cdots,0,z_n) \in \widetilde{\Omega},$ we have
    \begin{equation*}
         |v(0,\cdots,0,z_n)| \leq |W(0,\cdots,0,z_n)| = (n-1)(\frac{1}{\varepsilon_0})^{\theta} \cdot z_n^\theta .
    \end{equation*}
    In this way, we have proved \eqref{eq:theorem1-goal} and so Theorem \ref{theorem1}.

\section{The proof of Theorem 1.3}

With the assumption of Theorem 1.3,   we are going to prove  it.

    As the beginning of the proof of Theorem \ref{theorem1}, using $u(x)$ instead of $v(\widetilde{x})$ and $\Om$ instead of $\widetilde{\Omega}$
    we may assume that $x_0 = \boldsymbol{0} \in \partial\Omega$,  $\Omega$ satisfies
    $$ \{x\in\mathbb{R}^n \,| \, h> x_n > \eta_1 |x_1|^{a_1}+\cdots + \eta_{n-1} |x_{n-1}|^{a_{n-1}} \} \subseteq \Omega \subseteq \mathbb{R}_{+}^n, $$
    and  $ u\in  (\overline{\Omega})$ is a convex viscosity sub-solution to
    \begin{equation}
        \label{eq:theorem2-u}
        \left\{ \
        \begin{aligned}
        \det D^2 u & = |u|^{-n-k-2}((x+y_0)\cdot Du - u)^{-k}   \quad &in \quad & \Omega, \\
                u & = 0  \quad &on \quad &\partial \Omega,
      \end{aligned}
        \right.
    \end{equation}
    where $y_0 = (y_1,\cdots,y_n) \in \mathbb{R}^n $ such that $ |y_0| \leq \operatorname{diam} \Omega$. It is sufficient for Theorem 1.3 to prove
    \begin{equation}
        \label{neq:theorem2-goal}
        |u(0,\cdots,0,y_n)| \geq C\, y_n^{\,\theta}, \quad \forall y = (0,\cdots,0,y_n) \in \Omega,\ \ y_n\in(0,\frac{3}{8}h),
    \end{equation}
    where $ C = C(n,k,a_1,\cdots,a_{n-1}, \eta_1,\cdots,\eta_{n-1}, h,\operatorname{diam}\Omega) $.

    Here and below, we denote $x'= (x_1,\cdots,x_{n-1})$ for any $x\in \mathbb{R}^n$, and define the function
    $$v(x') = \eta(|x_1|^{a_1}+\cdots+|x_{n-1}|^{a_{n-1}})$$
     where $\eta = \max\{\eta_1,\cdots,\eta_{n-1}\}.$ Then $v(x')$ is smooth for $x' \neq 0.$ For any $\widetilde{h} \in (0,\frac{h}{2}),$ set
    \begin{equation}
        \begin{aligned}
            D' & = \{x' \in \mathbb{R}^{n-1} \mid v(x')< \widetilde{h}\} , \quad V=\{(x',x_n)\in \mathbb{R}^n \mid x'\in D', x_n = v(x') \}, \\
            V_1 & = \{(x',x_n)\in \mathbb{R}^n \mid v(x')< x_n < \widetilde{h} \}, \\
             V_2 &= \{(x',x_n)\in \mathbb{R}^n \mid 2v(x')< x_n < \widetilde{h} \},\\
            E & = \{(x^{\prime}, x_n) \in \mathbb{R}^n \mid \frac{x_1^2}{[c(\frac{\widetilde{h}}{\eta})^{\frac{1}{a_1}}]^2}+\cdots+\frac{x_{n-1}^2}{[c(\frac{\widetilde{h}}{\eta})^{\frac{1}{a_{n-1}}}]^2}+\frac{(x_n-\frac{3}{4} \widetilde{h})^2}{[c \widetilde{h}]^2} \leq 1\}.
        \end{aligned}
    \end{equation}
    By the convexity, we have $V \subset \Omega$, and $V_2 \subset V_1 \subset \Omega$. Furthermore, one can find a number $c(n)>0$ such that $E \subset V_2$ for all $c \in(0, c(n)]$. From now on, we fix $c=\min \left\{c(n), \frac{1}{4}\right\}$. Consequently
    \begin{equation}
        \label{neq:xn-h}
        \widetilde{h} \geq x_n \geq \frac{\widetilde{h}}{2}, \quad \forall x=\left(x^{\prime}, x_n\right) \in E .
    \end{equation}
    From the definition of $V$, there exists a constant $M:=C(n,a_1, \cdots, a_{n-1}, \eta, h) = C(\Omega)>0$ such that $\left|\nabla v\left(x^{\prime}\right)\right| \leq M$ for all $x^{\prime} \in D^{\prime}$ but $x^{\prime} \neq 0$, which together with the convexity implies
    $$
    d_x \geq \operatorname{dist}(x, V) \geq \frac{x_n-v\left(x^{\prime}\right)}{\sqrt{1+M^2}}, \quad \forall x=\left(x^{\prime}, x_n\right) \in E .
    $$
    Since $E \subset V_2$, then $v\left(x^{\prime}\right)<\frac{x_n}{2}$ for all $x=\left(x^{\prime}, x_n\right) \in E$, so we have
    $$
    d_x \geq \frac{x_n}{2 \sqrt{1+M^2}},\quad  \forall x=\left(x^{\prime}, x_n\right) \in E,
    $$
    which yields
    \begin{equation}
        \label{range-dx-xn}
        d_x \in [\frac{x_n}{2 \sqrt{1+M^2}} , \quad x_n ], \quad \forall x=\left(x^{\prime}, x_n\right) \in E.
    \end{equation}

    Now let $W \in C^2(E)\cap C(\bar{E})$ be the convex solution to the following problem
    \begin{equation} 
        \label{thm2-W-new}
        \begin{aligned}
            \det D^2 W & = C_1 \cdot \widetilde{h}^k \, |W|^{-n-2k-2}, \quad \forall x\in E, \\
            W|_{\partial E} & = 0.
        \end{aligned}
    \end{equation} where $C_1$ is a positive constant to be determined.
    The existence of $W$ can be obtained from \cite{[CY]}. Then we conclude that
     \begin{equation} 
     u\leq W\leq 0\ \ \text{ on}\ \  \bar E.
   \end{equation} Otherwise, by the assumption on $u$, we have some $\bar x\in E $ such that  $u-W$ attains a positive maximum at $\bar x$. Since $u$ is convex viscosity sub-solution
   to problem (4.1), by Definition 2.1 we have
    \begin{equation} 
    \begin{aligned}
   \det D^2W (\bar x) & \geq |u (\bar x)|^{-n-k-2}((\bar x+y_0) \cdot DW(\bar x) - u(\bar x))^{-k}\\
  & > | W(\bar x)|^{-n-k-2}((\bar x+y_0) \cdot DW(\bar x) - W(\bar x))^{-k} .
    \end{aligned}
    \end{equation}
    Since $u-W$ attains maximum at $\bar x$, we have
    $ u(\bar x)-W(\bar x)\geq u(x)-W(x)$ for all $x\in E.$
    Take any $p\in \partial u(\bar x)$. Then
    $$ u(x)\geq u(\bar x)+ p\cdot (x-\bar x)\geq W(\bar x)+u(x)-W(x)+ p\cdot (x-\bar x), \; \forall x\in E,$$
which yields
 $$ W(x)\geq   W(\bar x)+ p\cdot (x-\bar x), \; \forall x\in E.$$
 This means $p\in  \partial W(\bar x)=\{DW(\bar x)\}$, which implies  $\partial u(\bar x) =\{DW(\bar x)\}$. Therefore, we have
 $$ u(x)\geq   u(\bar x)+ DW(\bar x)\cdot (x-\bar x), \; \forall x\in E,$$ and the convexity of $u$ implies that
  $$ u(x)\geq   u(\bar x)+ DW(\bar x)\cdot (x-\bar x), \; \forall x\in \Om .$$
  Now we choose a unit  vector $\overrightarrow{e}$ such that $|D W(\bar x)|=D_{\overrightarrow{e}} W(\bar x)$, and  take a $z=x+|z-x| \overrightarrow{e} \in  \partial \Omega$.
  It follows from the above inequality that
  $$ DW (\bar x)\cdot (z-\bar x)\leq u(z)-u(\bar x)=-u(\bar x).$$ Combing (4.4) and (4.6), we have
    $$
    \begin{aligned}
    |D W(\bar x)| & =\left|D_{\vec{e}} W(\bar x)\right| \leq \frac{|u(\bar x)|}{|z-\bar x|} \leq \frac{|W(\bar x)|}{d_{\bar x}} \\
    & \leq \frac{|W(\bar x)|}{\frac{1}{2 \sqrt{1+M^2}} \cdot \bar x_n} \leq \frac{|W(\bar x)|}{\frac{1}{2 \sqrt{1+M^2}} \cdot \frac{\widetilde{h}}{2}} \\
    & =4 \sqrt{1+M^2} \cdot \frac{|W(\bar x)|}{\widetilde{h}}.
    \end{aligned}
    $$  This inequality, together with
      $|\bar x+y_0| \leq 2 \operatorname{diam}\Omega$ and  $\widetilde{h} < \frac{1}{2}h < \operatorname{diam}\Omega $,  implies
    \begin{equation*}
        \label{thm2-xDu-u}
        \begin{aligned}
            (\bar x+y_0)\cdot DW(\bar x) - W(\bar x) & \leq 2\operatorname{diam}\Omega \cdot 4 \sqrt{1+M^2} \cdot \frac{|W(\bar x)|}{\widetilde{h}} + |W(\bar x)| \\
            & \leq \frac{9\operatorname{diam}\Omega \sqrt{1+M^2}}{\widetilde{h}} |W(\bar x)| .
        \end{aligned}
    \end{equation*}
    Therefore, by this and (4.8) we obtain
    \begin{equation} 
        \begin{aligned}
         \det D^2W (\bar x) & > |W (\bar x)|^{-n-k-2}((\bar x+y_0) \cdot DW(\bar x) - W(\bar x))^{-k}\\
             & \geq |W(\bar x)|^{-n-k-2}[\frac{9\operatorname{diam}\Omega \sqrt{1+M^2}}{\widetilde{h}}]^{-k} |W(\bar x)|^{-k} \\
            & = C_1 \cdot \widetilde{h}^k |W(\bar x)|^{-n-2k-2},
        \end{aligned}
    \end{equation}
    where $C_1 = C(n,k,a_1,\cdots,a_{n-1},\eta,h,\operatorname{diam} \Omega).$  Obviously, (4.9) contradicts (4.6). In this way we have proved (4.7).
    As a special case of (4.7) we have
     \begin{equation}
        \label{u-W-h-tilde}
        |u(0,\cdots,0,\frac{3}{4}\widetilde{h})| \geq |W(0,\cdots,0,\frac{3}{4}\widetilde{h})| ,
    \end{equation}
    which will yield (4.2).  In fact, by  the   scaling transform
    $$
    z_1=\frac{x_1}{c(\frac{\widetilde{h}}{\eta})^{\frac{1}{a_1}}},\ \cdots, z_{n-1}=\frac{x_{n-1}}{c(\frac{\widetilde{h}}{\eta})^{\frac{1}{a_{n-1}}}},\ z_n=\frac{x_n-\frac{3}{4} \widetilde{h}}{c \widetilde{h}},\ \widetilde{W}(z)=\frac{W(x)}{(\widetilde{h})^\theta},
    $$
    where $\theta$ is given by \eqref{eq:theta-def1},  we see that
    \begin{equation*}
        \begin{aligned}
            \det D_z^2 \widetilde{W} & = (\widetilde{h}^{-\theta})^n \cdot \widetilde{h}^{\frac{2}{a_1}+\cdots+\frac{2}{a_{n-1}}+2}\, c^{2n}\, \eta^{-(\frac{2}{a_1}+\cdots+\frac{2}{a_{n-1}})} \det D_x^2 W \\
            & = \widetilde{h}^{-n\theta+\frac{2}{a_1}+\cdots+\frac{2}{a_{n-1}}+2} \, c^{2n}\, \eta^{-(\frac{2}{a_1}+\cdots+\frac{2}{a_{n-1}})}\, C_1 \, \widetilde{h}^k \, |W|^{-n-2k-2} \\
            & = C_1\, c^{2n}\, \eta^{-(\frac{2}{a_1}+\cdots+\frac{2}{a_{n-1}})}\, \widetilde{h}^{-n\theta+\frac{2}{a_1}+\cdots+\frac{2}{a_{n-1}}+2+k} \cdot \widetilde{h}^{(-n-2k-2)\theta}\, |\widetilde{W}|^{-n-2k-2} \\
            & = C_1\, c^{2n}\, \eta^{-(\frac{2}{a_1}+\cdots+\frac{2}{a_{n-1}})}\, \widetilde{h}^{ \frac{2}{a_1}+\cdots+\frac{2}{a_{n-1}}+2+k - (2n+2k+2)\theta} \, |\widetilde{W}|^{-n-2k-2}.
        \end{aligned}
    \end{equation*} It follows from (4.6) and
    \eqref{eq:theta-def1} that
    $$
    \begin{aligned}
    \operatorname{det} D_z^2 \widetilde{W} & = C_1\, c^{2n}\, \eta^{-(\frac{2}{a_1}+\cdots+\frac{2}{a_{n-1}})} \, |\widetilde{W}|^{-n-2k-2}, \quad \forall z \in B_1(\mathbf{0}), \\
    \widetilde{W}|_{\partial B_1(\mathbf{0})} & =0 .
    \end{aligned}
    $$
    Then we have
    $$
    \widetilde{W}(\mathbf{0}) = C(n,k,C_1, c, a_1, \cdots, a_{n-1},\eta) = C(n,k,a_1, \cdots, a_{n-1},\eta,\operatorname{diam}\Omega) <0 .
    $$
    Combining \eqref{u-W-h-tilde}, we obtain
    $$
    |u(0, \cdots, 0, \frac{3}{4} \widetilde{h})| \geq |W (0, \cdots, 0, \frac{3}{4} \widetilde{h} ) |= |(\widetilde{h})^\theta\, \widetilde{W}(\mathbf{0}) |= C \, (\frac{3}{4} \widetilde{h} )^\theta
    $$
    for all $\widetilde{h} \in (0, \frac{h}{2} )$ and $C$ is a positive constant depending only on $n,k,\operatorname{diam}\Omega$ and  $a_1, \cdots, a_{n-1},$ $\eta_1, \cdots, \eta_{n-1}, h$ in the interior convex condition at $x_0$. That is
    $$
    \left|u\left(0, \cdots, 0, y_n\right)\right| \geq C \, y_n^{\,\theta}, \quad \forall y_n \in(0, \frac{3}{8}h).
    $$
    In this way, we have proved \eqref{neq:theorem2-goal} and so Theorem \ref{theorem2}.

\section{The proof of Theorem 1.4}

 With the assumption of Theorem 1.4, we first prove its conclusion (1).  In this case, $u\in C(\overline{\Omega})$ is convex in  $\Omega$, and by Definitions 2.3 and 2.4 we have
   \begin{equation}\label{5.1}
   \begin{aligned}
   \mathbb{M} u & \leq  |u(x)|^qdx  \ \ \text{  in}\ \  \Om, \\
   u &=0\ \ \  \text{  on}\ \ \ \Om .
   \end{aligned}
   \end{equation}
   Following the proof of  Lemma 3.1 (i) in \cite{[Le1]} we see that
    \begin{equation}\label{5.2}
    ||u||_{L^\infty (\Om)}\leq \hat C(n, q)|\Om|^{\frac{2}{n-q}} .
    \end{equation}
       Since $ x_0 \in \partial\Omega$ is a  $(a_1, \cdots, a_{n-1})$-exterior  condition point and   \eqref{eq:monge2-q} is invariant under rotation and translation, we can assume that $x_0=\boldsymbol{0} \in \partial \Omega$ and
    $$ \Omega \subseteq \{x\in\mathbb{R}^n \,| \, x_n \geq \eta_1 |x_1|^{a_1}+\cdots+\eta_{n-1} |x_{n-1}|^{a_{n-1}} \}. $$
      It suffices to verify that for all $x\in \Omega$,
    \begin{equation} 
        \label{eq:xn-beta}
        |u(x)| \leq C(n,q,a_1,\cdots,a_{n-1},\eta_1,\cdots,\eta_{n-1},\operatorname{diam} \Omega,\lambda)\, x_n^\lambda,
    \end{equation} where $\lambda$ is the same as in Theorem 1.4 (1).

    \textbf{Case 1:} $n-\bar a-2>0$. In this case, $\lambda=\lambda_0= \frac{\bar a+2}{n}$. So $\lambda_0\in (0, 1)$.  Analogous to the proof of Theorem \ref{theorem1}, we consider
    \begin{equation}
        W(x)=- \sum_{i=1}^{n-1} \left[\left(\frac{x_n}{\varepsilon}\right)^{\frac{2}{a_i}}-x_i^2\right]^{\frac{1}{b_i}}, \ \ b_i = \frac{2}{a_i\lambda_0} , \quad \forall x \in \Omega,
        \label{eq:def-Wi-W-2}
    \end{equation}
    where $\varepsilon$ is a positive constant to be determined. Notice that the estimates \eqref{Wxixi} - \eqref{est-detD2W} still make sense. The only difference is $\lambda_0$ instead of $\theta$.  Therefore, we have
    \begin{equation}
        \label{eq:HWfinal2-2}
        \det D^2 W  \geq (\frac{1}{\varepsilon})^2 \cdot \tau_1(\varepsilon) \cdot(\frac{x_n}{\varepsilon})^{n \lambda_0-(\frac{2}{a_1}+\cdots+\frac{2}{a_{n-1}})-2} = (\frac{1}{\varepsilon})^2 \, \tau_1(\varepsilon),
    \end{equation}
    where $\tau_1(\varepsilon)$ satisfies
    $$
    \lim _{\varepsilon \to 0} \tau_1(\varepsilon) = (n-1)a_1 \cdots a_{n-1} \lambda_0 ^n(1-\lambda_0) >0 \ \ \text{ uniformly for } x \in \Omega.
    $$ Here we have used the fact $\lambda_0<1$.
    By this, (5.1), (5.2) and (5.5),  we can choose $\varepsilon = c(n,q,a_1,\cdots,a_{n-1},\eta_1,\cdots,\eta_{n-1},|\Om|)$ sufficiently small such that
    \begin{equation}
        \det D^2 W dx\geq  [\hat C(n, q)|\Om|^{\frac{2}{n-q}}]^q dx\geq \mathbb{M} u.
    \end{equation}
    Notice that $u=0\geq W$ on $\partial\Omega$, we have $u\geq W$ in $\Omega$ by Lemma 2.1. Therefore $|u| \leq |W|$, which  yields (5.3), since $|\Om|$ can be estimated by $\operatorname{diam}\Omega$.

\textbf{Case 2:}   $n-\bar a -2\leq 0$. In this case, $\lambda$ can be taken any number in $(0, 1)$.  Since $\Om$ is bounded, it is enough to prove (5.3) for
$\lambda \in (\frac{\bar a}{n}, 1)$.   Take $b_i = \frac{2}{a_i\lambda}$ for all $i \in \{1,\cdots,n-1\}$ in \eqref{eq:def-Wi-W-2}. As \eqref{eq:HWfinal2-2}, we get
    \begin{equation*}
        \det D^2 W(x) \geq (\frac{1}{\varepsilon})^2 \cdot \tau_1(\varepsilon) \cdot(\frac{x_n}{\varepsilon})^{n \lambda -\bar{a}-2} \geq (\frac{1}{\varepsilon})^2 \cdot \tau_1(\varepsilon) \cdot(\frac{\operatorname{diam}\Omega}{\varepsilon})^{n \lambda -\bar{a}-2},
    \end{equation*}
    where  $$
    \lim _{\varepsilon \to 0} \tau_1(\varepsilon) = (n-1)a_1 \cdots a_{n-1} \lambda^n (1-\lambda) >0 \ \ \text{ uniformly for } x \in \Omega.
    $$  Since $n \lambda -\bar{a}>0$ by the  assumption  $\lambda >\frac{\bar a}{n}$,
    choosing $\varepsilon$ sufficiently small and using the above estimate, (5.1) and (5.2),  we obtain
    $$   \det D^2 W >   [\hat C(n, q)|\Om|^{\frac{2}{n-q}}]^q dx\geq \mathbb{M} u.$$
     The rest of the proof is the same as Case 1.

     Now we are in the position to prove (2) of Theorem 1.4.  The assumption that  $u\in  C(\overline{\Omega})$
    is  Aleksandrov  solution to problem \eqref{eq:monge2-q} means  $u\in C^\infty (\Om)$. See Proposition 2.8 in \cite{[Le1]}.  Since $\Omega$ is $(a_1, \cdots, a_{n-1})$-exterior  convex domain,
     for any point $z \in \Omega$, let $x_0 \in \partial \Omega$ be a nearest boundary point to $z$. Since $\Omega$ is $(a_1, \cdots, a_{n-1})$-exterior  convex domain and \eqref{eq:monge2-q} is invariant under rotation and translation, we can assume that $x_0=\boldsymbol{0} \in \partial \Omega$ and
    $$ \Omega \subseteq \{x\in\mathbb{R}^n \,| \, x_n \geq \eta_1 |x_1|^{a_1}+\cdots+\eta_{n-1} |x_{n-1}|^{a_{n-1}} \}. $$
    This implies that $z$ lies on the $x_n$-axis. Denote $K := \| u \|_{L^\infty(\Omega)} $. By Lemma 3.1 (iii) in \cite{[Le1]}, there is a uniform estimate
    \begin{equation}
        \label{range-K-Omega}
        c(n,q)|\Omega|^{\frac{2}{n-q}} \leq K \leq C(n,q)|\Omega|^{\frac{2}{n-q}}.
    \end{equation}
    Using \eqref{range-K-Omega}, it suffices to verify that for all $x\in \Omega$,
    \begin{equation} 
        \label{eq:xn-beta-2}
        |u(x)| \leq C(n,q,a_1,\cdots,a_{n-1},\eta_1,\cdots,\eta_{n-1},K,\operatorname{diam} \Omega,\lambda)\, x_n^\lambda, \quad\quad \forall \lambda\in(0,\alpha).
    \end{equation}
\textbf{Case 1:} $ q\in (0, n-\bar a-2)$.  In this case $\alpha \in (0, 1)$.
In (1) we have proved \eqref{eq:xn-beta-2} for $\lambda=\lambda_0= \frac{\bar a+2}{n}$.  Now we use an iterative argument to prove it for each $\lambda \in (0, \alpha).$

 If for some $j \in \mathbb{N}, \lambda_j \in \left[\lambda_0, \alpha \right)$ and some positive constant $C_j$, we have
    \begin{equation}
        \label{eq:iterstart}
        |u(x)| \leq C_j \, x_n^{\lambda_j}, \quad\quad \forall x \in \Omega,
    \end{equation}
      then we can deduce that
    \begin{equation}
        \label{eq:iter}
        |u(x)| \leq C(n,q,a_1,\cdots,a_{n-1},\eta_1,\cdots,\eta_{n-1},\operatorname{diam} \Omega ,C_j ) \, x_n^{\lambda_{j+1}} ,  \quad\quad \forall x \in \Omega,
    \end{equation}
    where $\lambda_{j+1} := \frac{q \lambda_j}{n} + \frac{1}{n}(\bar{a} + 2) $. It is obvious that if $\lambda_j < \alpha $, then $\lambda_0 \leq \lambda_j < \lambda_{j+1} < \alpha $ and
    \begin{equation}
        \label{eq:jqn}
        \alpha - \lambda_{j+1} = \frac{q}{n}\left(\alpha - \lambda_j \right) .
    \end{equation}
    Take $$ W_{j+1} = - \sum_{i=1}^{n-1 }\left[ \left( \frac{x_n}{\varepsilon_{j+1}}\right)^{\frac{2}{a_i}}-x_i^2 \right]^{\frac{1}{b_{i,j+1}}}, \ \text{where }\ b_{i,j+1} = \frac{2}{a_i \lambda_{j+1}}.$$
    Following the procedures in (1), we can obtain
    \begin{equation}
        \label{eq:W_j+1}
            \det D^2 W_{j+1} \geq (\frac{1}{\varepsilon_{j+1}})^2 \cdot \tau_1(\varepsilon_{j+1}) \cdot(\frac{x_n}{\varepsilon_{j+1}})^{n \lambda_{j+1}-\bar{a}-2} \geq x_n^{n \lambda_{j+1}-\bar{a}-2}
    \end{equation}
    for some $\varepsilon_{j+1} $ sufficiently small.
    Suppose we have \eqref{eq:iterstart}, then we have
    \begin{equation}
        \label{eq:uC-hat}
        |u(x)|^{\frac{q}{n}} \leq (C_j)^{\frac{q}{n}} x_n^{\frac{q \lambda_j}{n}}< C' x_n^{\lambda_{j+1}-\frac{1}{n}(\bar{a}+2)}
    \end{equation}
    for some $C' = C(n,q,C_j) $ large enough.

    Let $U:=(\det D^2 u)(D^2 u)^{-1}$ be the cofactor matrix of the Hessian matrix $D^2 u$. Then
    $\det U=\left(\det D^2 u\right)^{n-1}$ and $\sum_{i, j=1}^n U^{i j} D_{i j} u=n \det D^2 u=n|u|^q$. Combining \eqref{eq:W_j+1} and \eqref{eq:uC-hat}, together with the matrix inequality
    $$ \operatorname{trace}(A B) \geq n(\det A)^{1 / n}(\det B)^{1 / n} \text{  for } A, B  \succeq  0, $$
    yields that
    \begin{equation}
        \label{eq:Uij-Dij}
        \begin{aligned}
           \sum_{i, j=1}^n U^{i j} D_{i j}\left(C' W_{j+1}\right) & \geq n C'\left(\operatorname{det} D^2 u\right)^{\frac{n-1}{n}}\left(\operatorname{det} D^2 W_{j+1}\right)^{\frac{1}{n}} \\
            & \geq n C' |u|^{\frac{q(n-1)}{n}} x_n^{\lambda_{j+1}- \frac{1}{n}(\bar{a}+2)} \\
            &>n|u|^q=n \operatorname{det} D^2 u = \sum_{i,j=1}^n U^{i j} D_{i j} u .
        \end{aligned}
    \end{equation}
    The maximum principle for the operator $\sum_{i, j=1}^nU^{i j} D_{i j}$ implies $u \geq C'\, W_{j+1}$ in $\Omega$, that is
    $$
    |u(x)| \leq |C' W_{j+1}(x)| \leq C_{j+1}\, x_n^{\lambda_{j+1}}, \quad\quad \forall  x \in \Omega ,
    $$
    where $ C_{j+1} = C(n,q,a_1,\cdots,a_{n-1},\eta_1,\cdots,\eta_{n-1},\operatorname{diam} \Omega ,C_j ). $ This gives \eqref{eq:iter}.

    Now we can choose $\lambda_0 = \frac{\bar{a}+2}{n}$ to initiate the iteration and obtain a sequence $\lambda_j$. From \eqref{eq:jqn}, we know that $\lambda_j$ is increasing and
    \begin{equation*}
        \alpha - \lambda_j = \left(\frac{q}{n}\right)^j \left(\alpha - \lambda_0\right) = \alpha \left(\frac{q}{n}\right)^{j+1} .
    \end{equation*}
    For any $\lambda \in (0,\alpha)$, we can choose $j$ large enough such that
    \begin{equation*}
        \alpha \left(\frac{q}{n}\right)^{j+1} < \alpha - \lambda.
    \end{equation*}
    Then we have $\lambda<\lambda_j<\alpha$ and for all $x\in\Omega$,
    \begin{equation*}
        |u(x)| \leq C_j\, x_n^{\lambda_j}\, \leq \, C(n,q,a_1,\cdots,a_{n-1},\eta_1,\cdots,\eta_{n-1},K,\operatorname{diam} \Omega,\lambda)\, x_n^\lambda,
    \end{equation*}
   which is   \eqref{eq:xn-beta-2} exactly.

      \textbf{Case 2:} $\max\{0, n-2-\bar a\}\leq q<n$.  By the result of (1), we need only to show the case  $0< n-2-\bar{a}\leq q<n$. We claim that {\sl
      if for some $\lambda\in [\frac{\bar{a}+2}{n},1)$
    \begin{equation}
        \label{eq:th2-claim-1}
        |u(x)|\leq C_{\lambda}\, x_n^\lambda  \quad\quad \forall x\in\Omega,
    \end{equation}
    then for any $ \lambda < \gamma < \min\{\frac{q}{n}\lambda + \frac{\bar{a}+2}{n},1\}$,
    \begin{equation}
        \label{eq:th2-claim-2}
        |u(x)|\leq C(n,q,a_1,\cdots,a_{n-1},\eta_1,\cdots,\eta_{n-1},\operatorname{diam} \Omega ,C_{\lambda} )\, x_n^\gamma,   \quad\quad \forall x\in\Omega.
    \end{equation}  }
    Notice that (1) shows that  \eqref{eq:th2-claim-1} holds for $\lambda=\frac{\bar{a}+2}{n}$.

    Suppose $\lambda\in [\frac{\bar{a}+2}{n},1)$ and $ \lambda < \gamma < \min\{\frac{q}{n}\lambda + \frac{\bar{a}+2}{n},1\}$, we have
    \begin{equation}
        \label{eq:beta-gamma}
        |u(x)|^{\frac{q}{n}} \leq (C_\lambda)^{\frac{q}{n}} \, x_n^{\frac{q}{n}\lambda} < C' \, x_n^{\gamma-\frac{\bar{a}+2}{n}},
    \end{equation}
    where $ C' = C(n,q,C_\lambda) $ sufficiently large. Set $b_i = \frac{2}{a_i \gamma}$ for all $i \in \{1,\cdots,n-1\}$ in \eqref{eq:def-Wi-W-2}, then \eqref{eq:HWfinal2-2} shows that $ \det D^2 W \geq x_n^{\,n\gamma - \bar{a}-2} $ for some $\varepsilon = \varepsilon(n,\gamma,a_1,\cdots,a_{n-1},\eta_1,\cdots,\eta_{n-1},\operatorname{diam} \Omega)$ sufficiently small. Similar to the analysis in \eqref{eq:Uij-Dij}, we have
    \begin{equation}
        \label{eq:Uij-Dij-2}
        \begin{aligned}
           \sum_{i, j=1}^n U^{i j} D_{i j}\left(C' W\right) & \geq n C' \left(\operatorname{det} D^2 u\right)^{\frac{n-1}{n}}\left(\operatorname{det} D^2 W\right)^{\frac{1}{n}} \\
            & \geq n C' |u|^{\frac{q(n-1)}{n}} x_n^{\gamma-\frac{\bar{a}+2}{n}} \\
            &>n|u|^q=n \operatorname{det} D^2 u =\sum_{i, j=1}^n U^{i j} D_{i j} u.
        \end{aligned}
    \end{equation}
    The maximum principle for the operator $\sum_{i, j=1}^nU^{i j} D_{i j}$ applied to $u$ and $ C' W$ gives $u \geq C' W $ in $\Omega$, that is
    $$
    |u(x)| \leq |C' W(x)| \leq C_\gamma \,x_n^{\gamma}, \quad\quad \forall  x \in \Omega,
    $$
    where $ C_{\gamma} = C(n,q,a_1,\cdots,a_{n-1},\eta_1,\cdots,\eta_{n-1},\operatorname{diam} \Omega ,C_\lambda ). $
    This gives \eqref{eq:th2-claim-2} and the claim is proved.

    From  Case 1 in Theorem 1.4 (1), we already have $u \in C^{0, \lambda_0}(\overline{\Omega})$, where $\lambda_0 = \frac{\bar{a}+2}{n} $. Let $ \lambda_{j+1} = \frac{q}{n}\lambda_j +\frac{\bar{a}+2}{n}, j \in \mathbb{N}$. Simple computation shows that
    \begin{equation*}
        \alpha - \lambda_j = \left(\frac{q}{n}\right)^j \left(\alpha - \lambda_0\right) = \alpha \left(\frac{q}{n}\right)^{j+1} .
    \end{equation*}
    Since $ q \geq n-2-\bar{a} $, we have $ \alpha \geq 1 $. Therefore, take any $\lambda\in(0,1)$, we can find a positive integer $j$ such that
    \begin{equation*}
        \alpha \left(\frac{q}{n}\right)^{j+1} < \alpha - \lambda.
    \end{equation*}
    Then we have $\lambda<\lambda_j$ and for all $x\in\Omega$,
    \begin{equation*}
        |u(x)| \leq C_{\lambda_j}\, x_n^{\lambda_j}\, \leq \, C(n,q,a_1,\cdots,a_{n-1},\eta_1,\cdots,\eta_{n-1},K,\operatorname{diam} \Omega,\lambda)\, x_n^\lambda,
    \end{equation*}
       which is   \eqref{eq:xn-beta-2} exactly. In this way, we complete the proof of Theorem 1.4.

\section{The proof of Theorem 1.5}

    As the beginning of the proof of (2) of Theorem 1.4,  we may assume that
    $$
    x_0=\mathbf{0} \text { and }\left\{\left.x \in \mathbb{R}^n\left|\eta_1\right| x_1\right|^{a_1}+\cdots+\eta_{n-1}\left|x_{n-1}\right|^{a_{n-1}}<x_n<h\right\} \subseteq \Omega \subseteq \mathbb{R}_{+}^n
    $$
    for some positive constants $\eta_1, \cdots, \eta_{n-1}$ and $h$.  Observing that in this case  $u\in  C(\overline{\Omega})\cap C^\infty (\Om)$
    is   smooth solution to problem \eqref{eq:monge2-q}  and $K := \Vert u \Vert_{L^{\infty}(\Omega)}$ satisfies  \eqref{range-K-Omega}. Hence,  we only need to prove for any $\lambda \in (\alpha,1)$,
    \begin{equation}\
        \label{eq:thm4-goal}
        \left|u\left(0,\cdots, 0, y_n\right)\right| \geq C \, y_n^\lambda, \quad \forall y=\left(0, \cdots, 0, y_n\right) \in \Omega,\ y_n \in(0, \frac{3}{8}h),
    \end{equation}
    where $C=C\left(n,q,a_1, \cdots, a_{n-1} , \eta_1, \cdots, \eta_{n-1} , h ,K, \operatorname{diam}\Omega ,\lambda \right).$

    \textbf{Step 1:} We will first prove estimate (6.1) for $\lambda_0:= \frac{1}{n}(q+\bar{a}+ 2)$. Assume $u(x_0) = -K$ where $x_0\in\Omega $. Then for any $ x\in\Omega$, the ray from $x_0$ to $x$ intersects $\partial\Omega$ at $z$. Write $x = \mu x_0 + (1-\mu)z $ for some $\mu \in (0,1]$. Since $u$ is a convex function and $u(z) = 0$, we have
    $$ u(x)\leq \mu u(x_0) + (1-\mu)u(z) = -\mu K, $$
    and therefore
    \begin{equation}
        \label{est-u-geq-dx}
        |u(x)| \geq \mu K = \frac{|x-z|}{|x_0-z|} K \geq \frac{\operatorname{dist}(x,\partial\Omega)}{\operatorname{diam}\Omega} K = \frac{K}{\operatorname{diam}\Omega} \, d_x.
    \end{equation}
    As in the proof of Theorem \ref{theorem2}, for any $\widetilde{h} \in (0,\frac{h}{2}),$ set
    \begin{equation}
        E = \{(x^{\prime}, x_n) \in \mathbb{R}^n \mid \frac{x_1^2}{[c(\frac{\widetilde{h}}{\eta})^{\frac{1}{a_1}}]^2}+\cdots+\frac{x_{n-1}^2}{[c(\frac{\widetilde{h}}{\eta})^{\frac{1}{a_{n-1}}}]^2}+\frac{(x_n-\frac{3}{4} \widetilde{h})^2}{[c \widetilde{h}]^2} \leq 1\}
    \end{equation}
    for some fixed $c = c(n)\leq \frac{1}{4}.$ Combining \eqref{neq:xn-h}-\eqref{range-dx-xn}, we get
    \begin{equation}
        \label{range-tildeh-dx}
         c_1\, \widetilde{h}\leq d_x\leq \widetilde{h}, \quad \forall x = (x',x_n) \in E,
    \end{equation}
    where $c_1 = C(n,a_1,\cdots,a_{n-1},\eta,h)>0. $ Consider
    \begin{equation}
        W(x) = - \varepsilon \left[1-\frac{x_1^2}{[c(\frac{\widetilde{h}}{\eta})^{\frac{1}{a_1}}]^2} - \cdots - \frac{x_{n-1}^2}{[c(\frac{\widetilde{h}}{\eta})^{\frac{1}{a_{n-1}}}]^2}- \frac{(x_n-\frac{3}{4} \widetilde{h})^2}{[c \widetilde{h}]^2} \right], \quad \forall x\in E,
    \end{equation}
    where $\varepsilon$ is a positive constant to be determined. Then $W|_{\partial E} = 0$ and
    \begin{equation}
        \det D^2 W = \frac{\varepsilon^n \cdot 2^n}{c^{2n}}\cdot \eta^{\bar{a}}\cdot \widetilde{h}^{-\bar{a}-2}.
    \end{equation}
    Noticing from \eqref{est-u-geq-dx} and \eqref{range-tildeh-dx}, we have
    \begin{equation}
        \det D^2 u = |u|^q \geq C(q,K,\operatorname{diam}\Omega)\,d_x^q \geq C(n,q,a_1,\cdots,a_{n-1},\eta,K,\operatorname{diam}\Omega) \, \widetilde{h}^q.
    \end{equation}
    Take $ \varepsilon = c_2 \, \widetilde{h}^{\frac{q+\bar{a}+2}{n}} = c_2 \, \widetilde{h}^{\lambda_0} $, where $c_2 = c(n,q,a_1,\cdots,a_{n-1},\eta,h,K,\operatorname{diam}\Omega)$ sufficiently small such that
    \begin{equation}
        \det D^2 u \geq C \, \widetilde{h}^q \geq \frac{\varepsilon^n \cdot 2^n}{c^{2n}}\cdot \eta^{\bar{a}}\cdot \widetilde{h}^{-\bar{a}-2} = \det D^2 W.
    \end{equation}
    Then maximum principle implies that $0  \geq W  \geq u$, which gives us
    $$|u(0, \cdots, 0, \frac{3}{4} \widetilde{h})| \geq |W (0, \cdots, 0, \frac{3}{4} \widetilde{h} ) |= \varepsilon= C_0 \cdot (\frac{3}{4} \widetilde{h} )^{\lambda_0}
    $$
    for all $\widetilde{h} \in (0, \frac{h}{2} )$ and $C_0$ is a positive constant depending only on $n,q, K,\operatorname{diam}\Omega$ and  $a_1, \cdots, a_{n-1},$ $\eta_1, \cdots, \eta_{n-1},h $ in the interior convex condition at $x_0$. That is
    $$
    \left|u\left(0, \cdots, 0, y_n\right)\right| \geq C_0 \cdot y_n^{\,\lambda_0}, \quad \forall y_n \in(0, \frac{3}{8}h).
    $$

    \textbf{Step 2:}  We claim that {\sl if for some $\lambda_j \in [ \alpha , \lambda_0 )$,
    \begin{equation}
        \label{eq:iterstart-2}
        |u(0,\cdots, 0, y_n)| \geq C_{\lambda_j} \, y_n^{\,\lambda_j} \quad\quad \forall y_n\in (0, \frac{3}{8}h),
    \end{equation}
    then
    \begin{equation}
        \label{eq:iter-2}
        |u(0,\cdots, 0, y_n)| \geq C_{\lambda_{j+1}} \, y_n^{\,\lambda_{j+1}} \quad\quad \forall y_n\in (0, \frac{3}{8}h),
    \end{equation}
    where $\lambda_{j+1} := \frac{q \lambda_j}{n} + \frac{1}{n}(\bar{a} + 2) $ and $ C_{\lambda_{j+1}} = C(n,q,a_1,\cdots,a_{n-1},\eta ,\operatorname{diam} \Omega ,C_{\lambda_j})$. }
    It is obvious that if $\lambda_j > \alpha $, then $\alpha  < \lambda_{j+1} < \lambda_{j} < \lambda_0 $ and
    \begin{equation}
        \label{eq:jqn-2}
        \lambda_{j+1} - \alpha = \frac{q}{n}\left(\lambda_j - \alpha \right) .
    \end{equation}
    Consider
    \begin{equation}
        W_{j+1}(x) = - \varepsilon_{j+1} \left[1-\frac{x_1^2}{[c(\frac{\widetilde{h}}{\eta})^{\frac{1}{a_1}}]^2}- \cdots -\frac{x_{n-1}^2}{[c(\frac{\widetilde{h}}{\eta})^{\frac{1}{a_{n-1}}}]^2}-\frac{(x_n-\frac{3}{4} \widetilde{h})^2}{[c \widetilde{h}]^2} \right], \quad \forall x\in E,
    \end{equation}
    where $\varepsilon_{j+1} $ is a positive constant to be determined. Following the procedures in \textbf{Step 1}, we can take $ \varepsilon_{j+1} = c_3 \, \widetilde{h}^{\frac{q \lambda_j+\bar{a}+2}{n}} = c_3 \, \widetilde{h}^{\lambda_{j+1}} $, where $ c_3 = c(n,q,a_1,\cdots,a_{n-1},\eta, \operatorname{diam}\Omega,C_{\lambda_j})$ sufficiently small such that
    \begin{equation}
        \det D^2 u = |u|^q \geq (C_{\lambda_j})^q \, \widetilde{h}^{q\lambda_j} \geq \frac{\varepsilon_{j+1}^n \cdot 2^n}{c^{2n}}\cdot \eta^{\bar{a}}\cdot \widetilde{h}^{-\bar{a}-2} = \det D^2 W_{j+1}.
    \end{equation}
    Then maximum principle implies that $0  \geq W_{j+1}\geq u $, which gives us
    $$|u(0, \cdots, 0, \frac{3}{4} \widetilde{h})| \geq |W_{j+1} (0, \cdots, 0, \frac{3}{4} \widetilde{h} ) |= \varepsilon_{j+1}= C_{\lambda_{j+1}} \cdot (\frac{3}{4} \widetilde{h} )^{\lambda_{j+1}},
    $$
    where $ C_{\lambda_{j+1}} = C(n,q,a_1,\cdots,a_{n-1},\eta ,\operatorname{diam} \Omega ,C_{\lambda_j} ). $ From the arbitrariness of $\widetilde{h} \in (0, \frac{h}{2})$, we have
    $$
    \left|u\left(0, \cdots, 0, y_n\right)\right| \geq C_{\lambda_{j+1}} \cdot y_n^{\,\lambda_{j+1}}, \quad \forall y_n \in(0, \frac{3}{8}h),
    $$
    which implies that \eqref{eq:iter-2} holds.

    Now we can choose $\lambda_0 = \frac{q+\bar{a}+2}{n}$ to initiate the iteration and obtain a suquence $\lambda_j$. From \eqref{eq:jqn-2}, we know that $\lambda_j$ is decreasing and
    \begin{equation*}
        \lambda_j - \alpha = \left(\frac{q}{n}\right)^j \left(\lambda_0 - \alpha\right) = (\alpha + 1) \left(\frac{q}{n}\right)^{j+1} .
    \end{equation*}
    For any $\lambda \in (\alpha,1)$, we can choose $j$ large enough such that
    \begin{equation*}
        (\alpha+1) \left(\frac{q}{n}\right)^{j+1} < \lambda - \alpha .
    \end{equation*}
    Then we have $\alpha<\lambda_j<\lambda$ and for all  $ y_n\in (0, \frac{3}{8}h),$
    \begin{equation*}
    |u(0,\cdots, 0, y_n)| \geq C_{\lambda_j} \, y_n^{\,\lambda_j}
         \geq  C(n,q,a_1, \cdots, a_{n-1} , \eta_1, \cdots, \eta_{n-1} , h ,K, \operatorname{diam}\Omega ,\lambda)\, x_n^{\lambda} .
    \end{equation*}
    Hence \eqref{eq:thm4-goal} holds and so is Theorem 1.5.

\end{document}